# Estimation of a Thermal Conductivity in a Stationary Heat Transfer Problem with a Solid-Solid Interface

Guillermo Federico Umbricht[1,2*], Diana Rubio[1], Domingo Alberto Tarzia[3]

[1] Centro de Matemática Aplicada, Escuela de Ciencia y Tecnología, Universidad Nacional de San Martín. 25 de Mayo y Francia, San Martín (B1650), Buenos Aires, Argentina

[2] Instituto de Ciencias e Instituto del Desarrollo Humano, Universidad Nacional de Gral. Sarmiento. Juan María Gutiérrez 1150, Los Polvorines (B1613), Buenos Aires, Argentina

[3] CONICET, Argentina - Departamento de Matemática, Facultad de Ciencias Empresariales, Universidad Austral, Paraguay 1950, S2000FZF Rosario, Santa Fe, Argentina

**ABSTRACT**

An inverse problem for a stationary heat transfer process is studied for a totally isolated bar on its lateral surface, made up of two consecutive sections of different, isotropic and homogeneous materials, perfectly assembly, where one of the materials, that is unreachable and unknown, has to be identified. The length of the bar is assumed to be much greater that the diameter so that a 1D heat transfer process is considered. A constant temperature is assumed at the end of the unknown part of the rod while the other end is let free for convection. We propose a procedure to identify the unknown material of the bar based on a noisy flow measurement at the opposite end. Necessary and sufficient conditions are derived together with a bound for the estimation error. Moreover, elasticity analysis is performed to study the influence of the data in the conductivity estimation and numerical examples are included to illustrate the proposed ideas and show the estimation performance.

## 1. INTRODUCTION

Thermal conductivity estimation problems under different conditions have been extensively studied in engineering for different applications. In some works, this parameter was determined from experimental, numerical and analytical techniques (see, e.g. [1],[2],[3],[4],[5],[6]). Other estimation approaches can be found in [7], [8], [9], [10].

On the other hand, heat transfer problems in multilayer or solid-solid interface materials have been studied in recent years due to the multiple and different applications in science and engineering. These problems have direct applications in several industries, including metallurgic, technologic, electronic, the automotive, aerospace and aviation [11], [12], [13], [14], [15], [16]. The advancement of technology requires materials with particular thermal, electrical, magnetic, acoustic and optical properties due to which the interface properties of different material combinations have been studied, for instance, Cu-Al [17], [18], Si-Ge [19], Al-Si [20], Pb-Sn [21], [22], Sn-Pb [23], Ti-Al [24].

Thermal conductivity is one of the most studied physical properties in heat transfer processes in composite materials as it can be seen in [25], [26], [27], [28], to mention some of the works found in the literature.

Under certain conditions, a process of heat transfer throughout a multilayer material can be seen as a process that develops in one direction through the different layers. In that case, it can be modeled as a 1D process with a solid-solid interface, see for instance, [29],[30].

In this work, we consider a stationary heat transfer process in a bar made up of two consecutive sections of different, isotropic and homogeneous materials, totally isolated on its lateral surface. The length of the bar is assumed to be much greater that the diameter so that a 1D heat transfer process is considered. The section on the left of the bar, which is assumed to be unreachable, is composed by an unknown material and its left end is assumed to have a constant temperature, represented by a Dirichlet condition. The thermal resistance at the interface is neglected and the right end of the bar is let free, in contact with an external fluid, giving rise to the phenomenon of convection, so a Robin type condition that models the heat dissipation is considered.

We propose to determine the thermal conductivity of unknown (left) material based on a noisy measurement of the heat flux at the available (right) body end. The approach presented here is quite simple and can be easily applied as a complement or alternative way for the determination of thermal conductivity in new materials.

Necessary and sufficient conditions for its application are derived together with an expression for a bound of the estimation error. In order to study the local influence of the heat flux measurement in the estimation, an elasticity analysis is performed and its properties are discussed. Some examples considering different situations are included to illustrate the estimation procedure proposed here. Absolute and relative errors are calculated for the results of the numerical experiments where noisy values of the heat flux are assumed for the estimation. In addition, we show that the errors obtained agree with the analysis of the elasticity function.

## 2. MATHEMATICAL MODEL

In this section, a mathematical model for the interface problem is stated and an analytical expression for the solution is given, which is consistent with the one corresponding to a homogenous bar. Numerical temperature profiles are shown for particular cases.

## 2.1 Model

The steady-state heat transfer problem for a bar composed of two consecutive, different, isotropic and homogeneous materials of known length $L$ (m) and diameter $d$ (m), $L \gg d$, fully insulated on its lateral surface, can be modeled by [30]

$$u''(x) = 0, \quad 0 < x < l, \quad (1)$$
$$u''(x) = 0, \quad l < x < L, \quad (2)$$

where $u$ (°C) represents the temperature of the bar and $l$ (m) the location of the contact point of the materials. A material A occupies the portion $(0, l)$ of the bar, and a material B occupies $(l, L)$.

Thermal resistance at the interface is assumed to be negligible; hence the following continuity conditions on temperature and heat flux are imposed [30]

$$u(l^+) = u(l^-), \quad (3)$$
$$\kappa_A u'(l^-) = \kappa_B u'(l^+), \quad (4)$$

where $\kappa_A, \kappa_B$ (Wm$^{-1}$°C$^{-1}$) represent the thermal conductivities of the materials A and B, respectively, $u(l^-) = \lim_{x \to l^-} u(x)$, $u(l^+) = \lim_{x \to l^+} u(x)$, and, analogously, $u'(l^-) = \lim_{x \to l^-} u'(x)$, $u'(l^+) = \lim_{x \to l^+} u'(x)$.

At the left boundary ($x=0$) is assumed to have a constant temperature, represented by a Dirichlet condition. The right boundary ($x=L$) is let free, in contact with an external fluid, giving rise to the phenomenon of convection, so a Robin type condition like a Newton law that models the heat dissipation is considered. Therefore, the boundary conditions are given by

$$u(x) = F, \quad x = 0, \quad (5)$$
$$\kappa_B u'(x) = -h(u(x) - T_a), \quad x = L, \quad (6)$$

where $F$ (°C) represents a constant temperature condition, $T_a$ (°C) the room temperature and $h$ (Wm$^{-2}$°C$^{-1}$) the convection heat transfer coefficient. In this work, it is assumed that $F > T_a$, analogous results may be obtained for $F < T_a$.

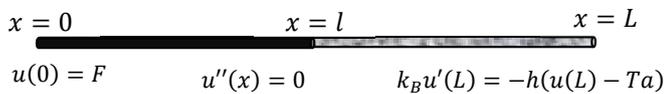

$x = 0$      $x = l$      $x = L$

$u(0) = F$    $u''(x) = 0$    $k_B u'(L) = -h(u(L) - Ta)$

**Figure 1.** Scheme for the mathematical model

## 2.2 Analytical solution of the Forward Problem

The solution to the problem described by Eqs. (1)-(6) is given by [31]

$$u(x) = F + \frac{\kappa_B h(T_a - F)}{\varsigma} x, \quad 0 \le x \le l, \quad (7)$$

$$u(x) = F + \frac{h(T_a - F)(l(\kappa_B - \kappa_A) + \kappa_A x)}{\varsigma}, \quad l < x \le L,$$

where

$$\varsigma = \kappa_A \kappa_B + \kappa_A h L + (\kappa_B - \kappa_A) h l. \quad (8)$$

Notice that if the bar is made of only one material, $\kappa_A = \kappa_B = \kappa$ then the solution (7)-(8) reduces to

$$u(x) = F + \frac{h(T_a - F)}{\kappa + hL} x, \quad 0 \le x \le L, \quad (9)$$

which is the solution for an homogenous bar [30].

## 2.3 Examples

Few examples are considered to illustrate the temperature profiles for the Eqs. (1)-(6) for different materials and contact points, by using the expressions (8)-(9). In all of them, it is assumed that $L = 1$ m, $F = 100$ °C, $T_a = 25$ °C and $h = 10$ Wm$^{-2}$°C$^{-1}$. Average thermal conductivity values were taken from [32] and included in Table 1.

Figure 2 shows the temperature profiles for a bar Fe-Cu and Cu-Fe with different interface locations. It can be seen that, when the interface is in the middle, the same temperature $u(L)$ is reached in both cases. In other words, if the bar is made up of equal parts of two materials, the location of these materials (left or right) does not influence the temperature value on the right end. This fact does not arise only from the numerical results but also from Eqs. (7) - (8) taking $l=L/2$, and it holds for any pair of materials.

**Table 1.** Average thermal conductivity values taken from [32] of the materials considered in this work

| Material | Symbol | $\kappa$ (Wm$^{-1}$°C$^{-1}$) |
|---|---|---|
| Aluminium | Al | 204 |
| Cupper | Cu | 386 |
| Iron | Fe | 73 |
| Silver | Ag | 419 |
| Lead | Pb | 35 |

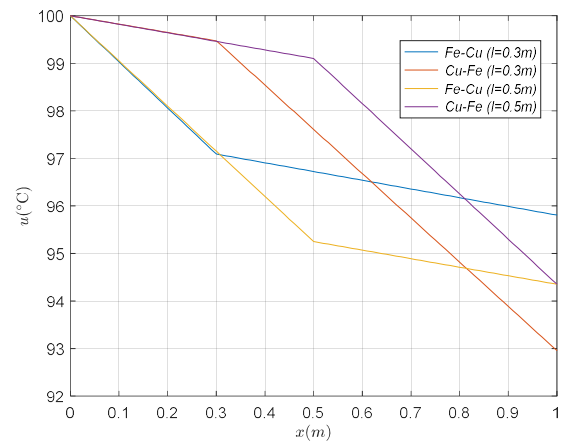

**Figure 2.** Temperature profiles for Fe-Cu and Cu-Fe at different solid-solid interface location.

Figure 3 shows the temperature profiles for three particular cases: $\kappa_A < \kappa_B$ (Fe-Ag), $\kappa_A > \kappa_B$ (Al-Pb) and $\kappa_A \approx \kappa_B$ (Ag-Cu) while keeping the same interface location. Note that the relationship between the thermal conductivity values $\kappa_A, \kappa_B$ is manifested in the angle that is formed in the temperature profile at the interface position. This property can also be obtained from Eq. (7) since

$$u(l^-) = u(l^+) = F + \frac{h(T_a - F)}{\varsigma} \kappa_B l,$$

then the angle at $x=l$ is given by

$$\alpha = \mathrm{accos}\left(\left(1 + \frac{a^2(\kappa_B - \kappa_A)^2}{(1 + a^2 \kappa_A \kappa_B)^2}\right)^{-1/2}\right)$$

where $a = \dfrac{h(T_a - F)}{\varsigma}$. Hence, $\alpha$ depends on $\kappa_B - \kappa_A$ and $\kappa_A \kappa_B$. Moreover, for $\kappa_A \approx \kappa_B$ it follows that $\alpha \approx \mathrm{accos}(1) = 0$.

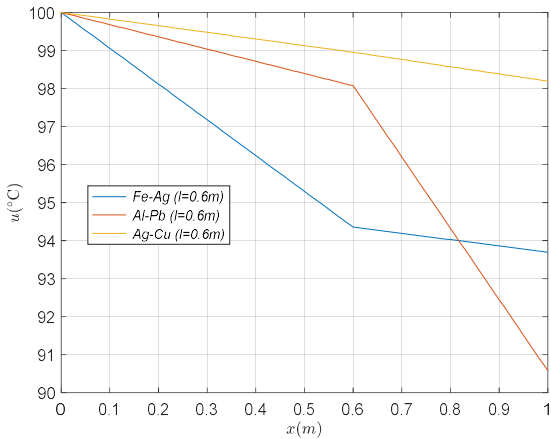

**Figure 3.** Temperature profiles for different materials with $\kappa_A < \kappa_B$ (Fe-Ag), $\kappa_A > \kappa_B$ (Al-Pb) and $\kappa_A \approx \kappa_B$ (Ag-Cu)

## 3. ESTIMATION OF THE UNKNOWN THERMAL CONDUCTIVITY

The main objective of this work is developed in this section. The thermal conductivity $\kappa_A$ of the material $A$ is estimated using a noisy heat flux data on $x = L$ subjected to Eqs. (1)-(6). Necessary and sufficient conditions for the estimation of the thermal conductivity are derived and a bound for the estimation error is provided.

### 3.1 Estimation

By definition, the thermal flux $q$ at $x = L$ is given by the following expression [29]

$$q = -\kappa_B u'(x), \quad x = L. \tag{10}$$

From Eq. (7), the thermal flux becomes

$$q = \frac{\kappa_B \kappa_A h(F - T_a)}{\kappa_A \kappa_B + \kappa_A h L + (\kappa_B - \kappa_A) h l}. \tag{11}$$

Then, the parameter $\kappa_A$ is obtained by the following expression,

$$\kappa_A = \frac{q h l \kappa_B}{h \kappa_B (F - T_a) - q h (L - l) - q \kappa_B}. \tag{12}$$

Therefore, the estimate $\hat{\kappa}_A$ of the thermal conductivity as a function of the heat flux measurement $\hat{q}$ at $x=L$ is expressed as

$$\hat{\kappa}_A = \frac{\hat{q} h l \kappa_B}{h \kappa_B (F - T_a) - \hat{q} h (L - l) - \hat{q} \kappa_B} \tag{13}$$

where it is assumed that

$$|q - \hat{q}| \leq \varepsilon, \tag{14}$$

and $\varepsilon > 0$ represents the noise level in the data.

Note that the estimation of $\hat{\kappa}_A$ depends on the heat flux measurement and the parameters of the problem.

### 3.2 Necessary and sufficient conditions

There exist necessary and sufficient conditions for the estimation of the thermal conductivity.

The thermal conductivity values $\kappa_A$, $\kappa_B$ and $\hat{\kappa}_A$ are assumed to be bounded, that is, there exist positive constants $\kappa_m$ and $\kappa_M$ that satisfy

$$0 < \kappa_m < \kappa_A, \kappa_B, \hat{\kappa}_A < \kappa_M \tag{15}$$

By Eq. (13) it follows that

$$0 < \kappa_m < \frac{\hat{q} h l \kappa_B}{h \kappa_B (F - T_a) - \hat{q} h (L - l) - \hat{q} \kappa_B} < \kappa_M. \tag{16}$$

Algebraic computations yield

$$q_m < q, \hat{q} < q_M < \bar{q}_M, \tag{17}$$

where

$$q_m = \frac{1}{l} \left( \frac{(F - T_a)}{\dfrac{1}{\kappa_m} + \dfrac{1 + \dfrac{h(L-l)}{\kappa_B}}{hl}} \right), \tag{18}$$

$$q_M = \frac{1}{l}\left(\frac{(F-T_a)}{\frac{1}{\kappa_M}+\frac{1+\frac{h(L-l)}{\kappa_B}}{hl}}\right) \quad (19)$$

and

$$\bar{q}_M = \frac{(F-T_a)h}{1+\frac{h(L-l)}{\kappa_B}}. \quad (20)$$

Therefore, there exists the estimate value $\hat{\kappa}_A$ and Eq. (15) is satisfied if and only if Eq. (17) holds. In other words, Eq. (17) is a necessary and sufficient condition for the existence of $\hat{\kappa}_A$.

### 3.3 Estimation error

Note that from Eqs. (12)-(13) it follows

$$\left|\frac{1}{\kappa_A}-\frac{1}{\hat{\kappa}_A}\right| = \frac{F-T_a}{lq\hat{q}}|q-\hat{q}|. \quad (21)$$

Eq. (17) leads to

$$\left|\frac{1}{\kappa_A}-\frac{1}{\hat{\kappa}_A}\right| < \left(\frac{F-T_a}{l}\right)\frac{|q-\hat{q}|}{q_m^2}, \quad (22)$$

which along with Eq. (18) it results

$$\left|\frac{1}{\kappa_A}-\frac{1}{\hat{\kappa}_A}\right| < \frac{|q-\hat{q}|\,l}{\kappa_m^2(F-T_a)}\left[1+\frac{\kappa_m}{hl}\left(1+\frac{h(L-l)}{\kappa_B}\right)\right]^2. \quad (23)$$

Equivalenty, by condition (14) we have

$$\left|\frac{1}{\kappa_A}-\frac{1}{\hat{\kappa}_A}\right| < K_1\varepsilon, \quad (24)$$

where the constant value $K_1$ is given by

$$K_1 = \frac{l}{\kappa_m^2(F-T_a)}\left[1+\frac{\kappa_m}{hl}\left(1+\frac{h(L-l)}{\kappa_B}\right)\right]^2. \quad (25)$$

Hence, $\left|\frac{1}{\kappa_A}-\frac{1}{\hat{\kappa}_A}\right| \to 0$ as $\varepsilon \to 0$.

Furthermore, an analytical bound can be determined for the estimation error of the coefficient of thermal conductivity. From (24)-(25) the estimation error $|\kappa_A - \hat{\kappa}_A|$ satisfies

$$|\kappa_A - \hat{\kappa}_A| < \kappa_A\hat{\kappa}_A K_1 |q-\hat{q}|. \quad (26)$$

or, by using the assumptions (14)-(15), the absolute estimation error satisfies

$$|\kappa_A - \hat{\kappa}_A| < K_2\varepsilon, \quad (27)$$

where the constant $K_2$ is given by the following expression

$$K_2 = \kappa_M^2 K_1 = \frac{\kappa_M^2 l}{\kappa_m^2(F-T_a)}\left[1+\frac{\kappa_m}{hl}\left(1+\frac{h(L-l)}{\kappa_B}\right)\right]^2. \quad (28)$$

Therefore, the estimation error $|\kappa_A - \hat{\kappa}_A|$ is bounded by $K_2\varepsilon$, and thus, $|\kappa_A - \hat{\kappa}_A| \to 0$ as $\varepsilon \to 0$.

## 4. ANALYSIS OF THE LOCAL DEPENDENCE OF THE PARAMETER WITH RESPECT TO THE DATA

Eq. (13) indicates that the estimated value $\hat{\kappa}_A$ for $\kappa_A$ depends on the parameters of the problem and on the measured heat flux $\hat{q}$. There are some tools that help to study the influence of data $\hat{q}$ on the estimated parameter $\hat{\kappa}_A$. Some of the most used are sensitivity [33] and elasticity analysis [34], depending on the discipline. In this work, the latter one is applied since it does not depend on the parameter to be estimated, as it is shown below.

### 4.1 Elasticity

This technique is widely used in economics. It provides the percentage error in the estimated parameter for 1 % error in a measurement value and it is defined by

$$E(q) = \frac{q}{\kappa_A}\frac{\partial \kappa_A}{\partial q}. \quad (29)$$

The expression (12) yields

$$E(q) = \frac{(F-T_a)h\kappa_B}{-q(\kappa_B+h(L-l))+h\kappa_B(F-T_a)}. \quad (30)$$

or, from (20),

$$E(q) = \frac{\bar{q}_M}{\bar{q}_M - q}. \quad (31)$$

### 4.2 Elasticity function analysis

The elasticity function given by the expression (31) does not depend on the parameter $\kappa_A$, and it has particularities that deserve to be highlighted.

#### 4.2.1 Vertical asymptote

The mathematical vertical asymptote for the function (30) (or (31)) is given by

$$q = \frac{h\kappa_B(F-T_a)}{h(L-l)+\kappa_B} = \bar{q}_M. \quad (32)$$

The vertical asymptote is located at $q = \bar{q}_M$. Although this value is never reached due to the restriction (17), the elasticity increases as $q$ approaches $\bar{q}_M$.

4.2.2 Positivity and monotonicity properties

Taking into account (17) and (31), we can ensure that the elasticity is a positive function for $q \in (q_m, q_M)$.

Another important observation is that the elasticity function is increasing. This fact can be easily seen by differentiating the expression (30) or (31) to obtain

$$\frac{\partial E(q)}{\partial q} = \frac{(F - T_a)h\kappa_B(\kappa_B + h(L-l))}{[-q(\kappa_B + h(L-l)) + h\kappa_B(F - T_a)]^2} > 0, \quad (33)$$

or, equivalently,

$$\frac{\partial E(q)}{\partial q} = \frac{\bar{q}_M}{(\bar{q}_M - q)^2} > 0. \quad (34)$$

## 5. NUMERICAL EXAMPLES

The following parameter values are imposed in the following numerical experiments $L = 10$ m, $l = 4$ m, $F = 100\,°C$, $T_a = 25\,°C$, $h = 10$ W.m$^{-2}$.°C$^{-1}$.

For the error testing, a few bars composed of different pairs of materials are considered where the thermal conductivity of the material occupying the left side must be estimated. The different pairs of materials are chosen according to the different relationships between the coefficients of thermal conductivity, i.e., $\kappa_A < \kappa_B$, $\kappa_A > \kappa_B$ and $\kappa_A \approx \kappa_B$.

In order to simulate a noisy experimental measurement $\hat{q}$, the solution $u$ to the forward problem (1)-(6) is calculated by using expressions (7)-(8) and the heat flux $q$ is obtained by (10). Finally, a perturbation is added to $q$ so that the resulting value remains within the interval given by condition (17).

The estimate value $\hat{\kappa}_A$ is calculated by using the expression (13) for the data $\hat{q}$. Afterwards, the absolute and relative estimation errors are calculated. This procedure is repeated for ten data perturbated with different noise levels.

For each example, a table with the results and a figure showing the elasticity function are included.

Each line in the tables below contain the noisy heat flux $\hat{q}$ along with the data error $|q - \hat{q}|$, the estimate $\hat{\kappa}_A$, the absolute error $|\kappa_A - \hat{\kappa}_A|$ and the relative error $\frac{|\kappa_A - \hat{\kappa}_A|}{\kappa_A}$.

### 5.1 Example 1 ($\kappa_A < \kappa_B$)

A Fe-Ag bar is considered to test the estimation procedure for the thermal conductivity of Fe, $\kappa_A = \kappa_{Fe} = 73$ W.m$^{-1}$.°C$^{-1}$. Since $\kappa_B = \kappa_{Ag} = 419$ W.m$^{-1}$.°C$^{-1}$, the analytical heat flux, calculated by using the expression (11), is $q = 443.48$ W.m$^{-2}$.

Table 2 shows the performance for the estimation procedure for this particular case. It is observed that all relative errors are less than 4% for noise level less 1%.

The elasticity function is plotted in Figure 4. As noted before, the elasticity function is positive and strictly increasing in the interval $(q_m, \bar{q}_M)$ where $\bar{q}_M = 656.05$ W.m$^{-2}$ is the vertical asymptote (31)-(32). This function indicates that a measurement error of 1 % in the heat flux value $\hat{q}$ translates into an error of around 4 % in the estimation value, which was also noticed before from the results shown in Table 2.

**Table 2.** Example 1: Estimates for $\kappa_A = 73$ W.m$^{-1}$.°C$^{-1}$, where the analytical heat flux is $q = 443.48$ W.m$^{-2}$.

| Data | Estimated value | Data error | Absolute estimation error | Relative estimation error |
|---|---|---|---|---|
| $\hat{q}$ | $\hat{\kappa}_A$ | $|q - \hat{q}|$ | $|\kappa_A - \hat{\kappa}_A|$ | $|\kappa_A - \hat{\kappa}_A|/\kappa_A$ |
| 439 | 70.767 | 4.480 | 2.232 | 0.031 |
| 440 | 71.257 | 3.480 | 1.742 | 0.024 |
| 441 | 71.751 | 2.480 | 1.248 | 0.017 |
| 442 | 72.249 | 1.480 | 0.750 | 0.010 |
| 443 | 72.753 | 0.480 | 0.246 | 0.003 |
| 444 | 73.261 | 0.520 | 0.261 | 0.004 |
| 445 | 73.774 | 1.520 | 0.774 | 0.011 |
| 446 | 74.292 | 2.520 | 1.292 | 0.018 |
| 447 | 74.814 | 3.520 | 1.815 | 0.025 |
| 448 | 75.342 | 4.520 | 2.342 | 0.032 |

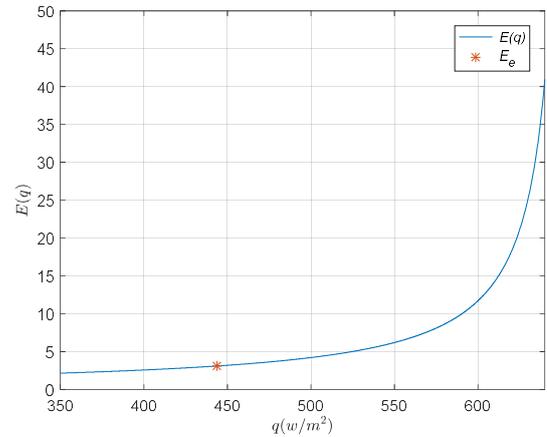

**Figure 4.** Elasticity function for Example 1.

### 5.2 Example 2 ($\kappa_A > \kappa_B$)

A Al-Pb bar is now considered for the estimation of the thermal conductivity value $\kappa_A = \kappa_{Al} = 204$ W.m$^{-1}$.°C$^{-1}$. For this example $\kappa_B = \kappa_{Pb} = 419$ W.m$^{-1}$.°C$^{-1}$ hence, the relationship between the thermal conductivities is $\kappa_A > \kappa_B$. For this example, the analytical heat flux is $q = 257.69$ W.m$^{-2}$. Ten simulated data are considered around the $q$ in the interval $(q_m, \bar{q}_M)$ and the results are shown in Table 3. The relative error for this example achieves 25% for an error of 2% in the data. The elasticity function for the Example 2 is plotted in the Figure 5 and it indicates that a measurement error of 1 % in the heat flux value $\hat{q}$ produces an error of 15 % in the estimation value $\hat{\kappa}_A$. Note that for this example the vertical asymptote is located at $\bar{q}_M = 276.31$ W.m$^{-2}$, not too far from the value of $q$, hence the estimation is not as good as for the previous example.

**Table 3.** Example 2: Estimates for $\kappa_A = 204$ W.m$^{-1}$.°C$^{-1}$, where the analytical heat flux is $q = 257.69$ W.m$^{-2}$.

| Data | Estimated value | Data error | Absolute estimation error | Relative estimation error |
|---|---|---|---|---|
| $\hat{q}$ | $\hat{\kappa}_A$ | $\|q - \hat{q}\|$ | $\|\kappa_A - \hat{\kappa}_A\|$ | $\|\kappa_A - \hat{\kappa}_A\|/\kappa_A$ |
| 252 | 152.728 | 5.690 | 51.272 | 0.251 |
| 253 | 159.909 | 4.690 | 44.090 | 0.216 |
| 254 | 167.296 | 3.690 | 32.264 | 0.158 |
| 255 | 176.296 | 2.690 | 27.703 | 0.136 |
| 256 | 185.699 | 1.690 | 18.300 | 0.090 |
| 257 | 196.076 | 0.690 | 7.923 | 0.039 |
| 258 | 207.586 | 0.310 | 3.586 | 0.018 |
| 259 | 220.425 | 1.310 | 16.425 | 0.081 |
| 260 | 234.838 | 2.310 | 30.838 | 0.151 |
| 261 | 251.134 | 3.310 | 47.134 | 0.231 |

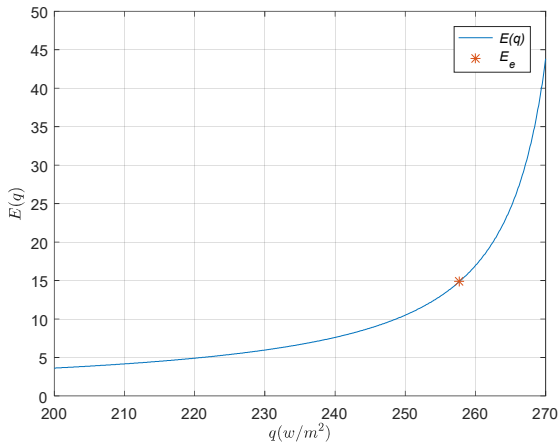

**Figure 5.** Elasticity function for Example 2.

### 5.3 Example 3 ($\kappa_A \approx \kappa_B$)

Finally, the case $\kappa_A \approx \kappa_B$ is considered. The thermal conductivity estimation is analyzed for a Ag-Cu bar is where $\kappa_A = \kappa_{Ag} = 419$ W.m$^{-1}$.°C$^{-1}$ and $\kappa_B = \kappa_{Cu} = 386$ W.m$^{-1}$.°C$^{-1}$.

The heat flux in this case is $q = 599.56$ W.m$^{-2}$ and from the necessary and sufficient condition (17), the upper bound for the measured data is $\bar{q}_M = 649.10$ W.m$^{-2}$. Ten simulated data are considered around the $q$, within the interval $(q_m, \bar{q}_M)$. The results are shown in Table 4. The relative error for the set of data considered here almost reaches 11%. The elasticity function for this example, shown in Figure 6, indicates that a measurement error of 1% in the heat flux value $\hat{q}$ produces an estimation error of around 13%. Again, as in the previous example, the value of $q$ is close to the vertical asymptote value $\bar{q}_M$ and therefore, the error in the data value is amplified in the estimation.

**Table 4.** Example 3: Estimates for $\kappa_A = 419$ W.m$^{-1}$.°C$^{-1}$, where the analytical heat flux is $q = 599.56$ W.m$^{-2}$.

| Data | Estimated value | Data error | Absolute estimation error | Relative estimation error |
|---|---|---|---|---|
| $\hat{q}$ | $\hat{\kappa}_A$ | $\|q - \hat{q}\|$ | $\|\kappa_A - \hat{\kappa}_A\|$ | $\|\kappa_A - \hat{\kappa}_A\|/\kappa_A$ |
| 595 | 380.721 | 4.690 | 38.278 | 0.091 |
| 596 | 388.542 | 3.690 | 30.457 | 0.073 |
| 597 | 396.664 | 2.690 | 22.336 | 0.053 |
| 598 | 405.103 | 1.690 | 13.896 | 0.033 |
| 599 | 413.879 | 0.690 | 5.120 | 0.012 |
| 600 | 423.013 | 0.310 | 4.014 | 0.010 |
| 601 | 432.527 | 1.310 | 13.527 | 0.032 |
| 602 | 442.527 | 2.310 | 23.445 | 0.056 |
| 603 | 452.792 | 3.310 | 33.792 | 0.081 |
| 604 | 463.599 | 4.310 | 44.599 | 0.106 |

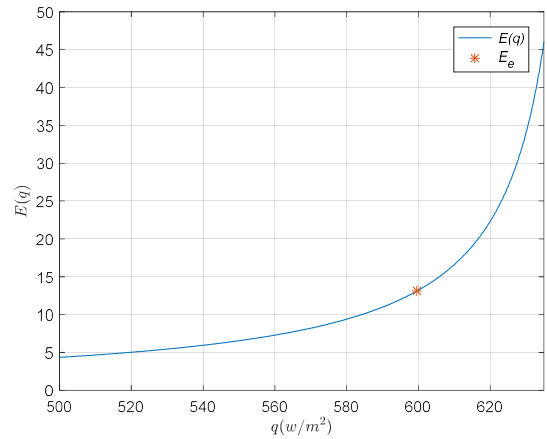

**Figure 6.** Elasticity function for Example 3

### 6. CONCLUSIONS

A mathematical model for the heat transfer along a bar totally isolated on its lateral surface composed of two different, isotropic and homogeneous materials, is considered. It is assumed that a unidimensional mathematical model can be used to describe the process. Appropriate boundary conditions are imposed and an analytical solution to the forward problem is found.

An procedure is proposed for the estimation of the thermal conductivity of the material that occupies the left part of the bar, using an over-specified heat flux condition at the right boundary. Necessary and sufficient conditions are provided as well as a bound for the estimation error. Furthermore, the local dependence of the data on the estimate is analyzed by means of the elasticity function, which provides a way to measure how the measurement error of the data influences the estimation error.

Numerical examples are included illustrating different set-ups, according to the relationship between the thermal conductivities of the materials, that is, $\kappa_A < \kappa_B$ (Fe-Ag), $\kappa_A > \kappa_B$ (Al-Pb) and $\kappa_A \approx \kappa_B$ (Ag-Cu), for the same boundary and initial conditions and the same bar length and interface solid-solid location. For these numerical experiments, the relative error in the data was amplified, obtaining a greater relative error in the estimation. It was also observed, that the elasticity function is positive and increasing, and it has a vertical asymptote which is an upper bound for the heat flux value. Hence, as the heat flux approaches that value, a worse estimation is obtained. As the elasticity function is independent of the parameter $\kappa_A$ to be estimated, it can be calculated and analyzed before performing the estimation. Therefore, it is convenient to study this function before carrying out the estimation to know if the estimated value will have the desired precision.

## ACKNOWLEDGMENT

The first and second authors acknowledge support from SOARD/AFOSR through grant FA9550-18-1-0523. The third author acknowledges support from European Union's Horizon 2020 Research and Innovation Programme under the Marie Sklodowska-Curie Grant Agreement No. 823731 CONMECH and by the Project PIP No. 0275 from CONICET-UA, Rosario, Argentina.

## NOMENCLATURE

| | |
|---|---|
| $E$ | Elasticity |
| $F$ | heat source, °C |
| $H$ | convective coefficient, W. m$^{-2}$.°C$^{-1}$ |
| $L$ | interface position, m |
| $L$ | bar length, m |
| $Q$ | thermal flow, W. m$^{-2}$ |
| $\hat{q}$ | measured thermal flux, W. m$^{-2}$ |
| $\bar{q}$ | asymptote thermal flux, W. m$^{-2}$ |
| $T_a$ | room temperature, °C |
| $U$ | bar temperature, °C |
| $X$ | special variable, m |

**Greek symbols**

| | |
|---|---|
| $\varepsilon$ | bound for flow measurement error, W. m$^{-2}$ |
| $\kappa$ | thermal conductivity, W.m$^{-1}$.°C$^{-1}$ |
| $\hat{\kappa}$ | estimated thermal conductivity, W.m$^{-1}$°C$^{-1}$ |
| $\varsigma$ | auxiliary assistant, W$^2$.m$^{-2}$.°C$^{-2}$ |

**Subscripts**

| | |
|---|---|
| $A$ | material A |
| $B$ | material B |
| $e$ | exact value |
| $m$ | Minimum |
| $M$ | maximum |